\documentclass[letterpaper, 10 pt, conference]{ieeeconf}

\IEEEoverridecommandlockouts         % This command is only
\usepackage{amsmath}
\usepackage{amssymb}
\usepackage{mathrsfs}
\usepackage{bbm}
\usepackage{mathtools}
\usepackage[noend]{algpseudocode}
\usepackage{algorithm}
\usepackage{graphicx}
\usepackage[utf8]{inputenc}
\usepackage{units}
\floatname{algorithm}{Procedure}

\usepackage{xcolor}
\usepackage{xifthen}

\renewcommand{\t}{^{\mbox{\tiny \sf T}}}
\newcommand{\addt}{\mbox{\tiny \sf T}}

\newcommand{\inv}{^{- 1}}

\DeclareMathOperator{\E}{E}
\DeclareMathOperator{\tr}{tr}

\DeclareMathOperator{\ncdf}{cdfn}
\renewcommand{\d}{\mathrm{d}}
\let\Pr\relax
\newcommand{\Pr}{\mathbb{P}}
\newcommand{\defeq}{\vcentcolon=}

\renewcommand{\Re}{\mathbb{R}}
\newcommand{\N}{\mathbb{N}}
\renewcommand{\E}{\mathbb{E}}

\newcommand{\Amat}{\mathcal{A}}
\newcommand{\Bmat}{\mathcal{B}}
\newcommand{\Gmat}{\mathcal{G}}

\newcommand{\Amati}{\mathcal{A}^i}
\newcommand{\Bmati}{\mathcal{B}^i}

\newcommand{\Amatit}{\mathcal{A} ^ {i \mbox{\tiny \sf T}} }

\newcommand{\norm}[1]{\left\lVert#1\right\rVert}

\renewcommand{\Re}{\mathbb{R}}

\DeclareSymbolFont{matha}{OML}{txmi}{m}{it}% txfonts
\DeclareMathSymbol{\varv}{\mathord}{matha}{118} % sharp v
\DeclareMathSymbol{\varw}{\mathord}{matha}{119} % sharp w

\newtheorem{problem}{Problem}
\newtheorem{remark}{Remark}
\graphicspath{ {./img/} }

\title{\LARGE \bf
Nonlinear Uncertainty Control with Iterative Covariance Steering
}

\author{Jack Ridderhof \qquad Kazuhide Okamoto \qquad Panagiotis Tsiotras%
\thanks{This work of the first author was supported by NASA Space Technology Research Fellowship 80NSSC17K0093.
The work of the second and third authors was supported by NSF award CPS-1544814.
The second author was also partially supported by the Funai Foundation for Information Technology.}% <-this % stops a space
\thanks{J. Ridderhof is a PhD student with the School of Aerospace Engineering, Georgia Institute of Technology, Atlanta, GA, 30332-0150, USA. Email: jridderhof3@gatech.edu}%
\thanks{K. Okamoto is a PhD student with the School of Aerospace Engineering, Georgia Institute of Technology, Atlanta, GA, 30332-0150, USA. Email: kazuhide@gatech.edu}%
\thanks{P. Tsiotras is the Andrew and Lewis Chair Professor with the D. Guggenheim School of Aerospace Engineering and the Institute for Robotics and Intelligent Machines, Georgia Institute of Technology, Atlanta, GA 30332-0150, USA. Email: tsiotras@gatech.edu}
}

\begin{document}

\maketitle
\thispagestyle{empty}
\pagestyle{empty}

%%%%%%%%%%%%%%%%%%%%%%%%%%%%%%%%%%%%%%%%%%%%%%%%%%%%%%%%%%%%%%%%%%%%%%%%%%%%%%%%
\begin{abstract}

%This paper considers the problem of steering the state distribution of a nonlinear stochastic system from an initial Gaussian to a specified terminal Gaussian, subject to probabilistic path constraints. An algorithm is developed to solve this problem by iteratively solving an approximate linearized problem, which we call the iterative covariance steering (iCS) method.

%This paper considers the problem of steering the state distribution of a nonlinear stochastic system from an initial Gaussian to a terminal distribution with a specified mean and covariance, subject to probabilistic path constraints. An algorithm is developed to solve this problem by iteratively solving an approximate linearized problem, which we call the iterative covariance steering (iCS) method.

This paper considers the problem of steering the state distribution of a nonlinear stochastic system from an initial Gaussian to a terminal distribution with a specified mean and covariance, subject to probabilistic path constraints. An algorithm is developed to solve this problem by iteratively solving an approximate linearized problem as a convex program. This method, which we call iterative covariance steering (iCS), is numerically demonstrated by controlling a double integrator with quadratic drag force subject to additive Brownian noise while satisfying probabilistic path constraints.

\end{abstract}

%%%%%%%%%%%%%%%%%%%%%%%%%%%%%%%%%%%%%%%%%%%%%%%%%%%%%%%%%%%%%%%%%%%%%%%%%%%%%%%%
\section{INTRODUCTION}

% - Coupling between feedback and trajectory (call out seperation principle for linear systems?)
% - State constraints included in coupling.
% - Model as stochastic optimal control problem, naturally have notion of coupling
% - Now BCs become distributions and constraints become chance constrints
% - Linear case - covariance steering
% - Nonlinear case - applications, proposed method
% - Contribution, summary of paper

%\textcolor{red}{Want to add: reference to successive convexifiation, DOC, more background and citations, and two motivating applications (powered descent, proximity operations with 3 body problem dynamics)}

Guidance and control design has generally followed the standard approach where an open-loop reference optimal control is solved with respect to the nonlinear dynamics and then a feedback controller is subsequently designed with respect to the dynamics linearized about the reference trajectory. Hence, there is an implicit unidirectional dependence of the feedback controller on this reference trajectory, but there is no direct mechanism from which the reference trajectory is affected by the closed-loop system behavior. Intuitively, if we could explicitly couple the design of the reference trajectory with the design of the feedback controller, then since we are optimizing over a larger set we could improve closed-loop system performance.
% Need to add more to this:
The situation becomes more complex with the introduction of state constraints and uncertainty. If the closed-loop statistics of a system are not considered, then the reference trajectory design must be conservative to satisfy the constraints. For systems that are significantly influenced by uncertain external forces, the conservatism of this approach may lead to greatly increased control cost or even infeasibility.

%One proposed method to capture this coupling between guidance and control is by a first order approximation of state perturbations from the reference trajectory by means of the sensitivity function. The optimal trajectory is then computed with an augmented cost function that includes a sensitivity function \cite{Seywald1996, Shen2010}.

In this paper we consider the system state to be a random vector which evolves according to a nonlinear stochastic differential equation with additive Brownian noise.
By letting the state to be a random vector, the control problem can be formulated as one of simultaneously steering each sample trajectory, and as a consequence, we can analytically study the difference between open and closed-loop control \cite{Brockett2010}. This machinery will serve as our mechanism to understand the coupling between the reference trajectory and the feedback controller.
%Furthermore, since the state is random, we consider boundary conditions on the state distribution. Specifically, we assume that the state is normally distributed and steer the mean and covariance to specified terminal values.
We assume that the state is normally distributed at the initial time, and we will design a control that steers the mean and the covariance of the initial state distribution to some given terminal values at the final time. This problem is referred to as the nonlinear covariance steering (CS) problem.
Since the state is assumed to be normally distributed, and thus unbounded, we must treat state constraints probabilistically. That is, the probability that the constraints are satisfied must be greater than some prespecified value. Since, by construction, these constraints may not be met for every sample path, they are often referred to as chance constraints \cite{Charnes1963, Nemirovski2006}.

%We are interested in guiding a spacecraft to a soft landing on the surface. In order to add robustness to random disturbances, this problem was posed as a stochastic optimal control problem \cite{Ridderhof2019}. However, current human Mars mission architectures include powered descent at supersonic speeds, and so the guidance design must take into account nonlinear aerodynamic drag terms \cite{Cianciolo2017}.

%Another problem receiving significant attention is spacecraft rendezvous and proximity operations. It has been presented as a stochastic optimal control problem. \cite{Mammarella2017} Plans to construct a cis-lunar space station (cite) have renewed interest in spacecraft rendezvous guidance with three-body-problem dynamics. \cite{Franzini2017}

The special case of linear time-varying stochastic systems with additive Brownian noise has been extensively studied in the literature. It has been shown that if the system is controllable, then the state covariance is also controllable \cite{Brockett2010}. That is, for an initial covariance $P_{x_0} > 0$ at time $t_0$, there exist a state feedback gain defined on the interval $[t_0, t_f]$ that steers the covariance to any final value $P_{x_f} > 0$ for any time $t_f > t_0$. The solution to the optimal linear CS problem with expected quadratic cost was given by Chen et al. \cite{Chen2016a, Chen2016b, Chen2016c}, and the solution was found to be closely related to the classical linear quadratic feedback control. The discrete linear CS problem with quadratic cost has also been studied and a similar close connection to linear quadratic control has been shown \cite{Goldshtein2017}.

It is well known that, for linear systems, and in the absence of any constraints, the mean and the covariance have independent dynamics, and that the mean state is controlled by the mean open-loop control and the covariance is controlled by the state feedback gain.
%Hence, we can consider the mean steering and the covariance steering as separate problems. 
It follows that, without constraints, we can consider the mean steering and the covariance steering as separate problems but that, when there are constraints, the mean and covariance are coupled.
%However, the mean and covariance are coupled through the constraints. 
In other words, for linear systems, the reference trajectory explicitly depends on the closed-loop behavior of the system when the state or control is constrained. 
For the discrete-time linear case with convex chance constraints, the chance-constrained CS problem can be cast as a deterministic convex optimization problem \cite{Okamoto2018b}.
This work was later extended to include non-convex chance constraints using mixed-integer programming~\cite{Okamoto2018}.

%%%

In this paper, we propose an algorithmic approach to solve the nonlinear CS problem by iteratively solving the linear CS problem with respect to the reference trajectory of the previous step.
%
%This algorithm, which we will refer to as iterative CS (iCS), is natural extension of linear CS in the spirit of other well known successive approximation methods such as differential dynamic programming (DDP) (ref) and iterative LQG (iLQG) (ref). While both DDP and iLQG compute a feedback control by backwards propagating an approximation of the value function, iCS solves each iteration as a constrained convex program. On one had, this has the benefit of more directly incorporating constraints, but comes at the cost of computation time and only approximating the dynamics to first order.
%
This algorithm, which we will refer to as iterative CS (iCS), is a natural extension of linear CS in the spirit of other well known successive approximation methods such as differential dynamic programming (DDP) \cite{Theodorou2010} and iterative LQG (iLQG) \cite{Todorov2005}, which both compute a feedback control by backwards propagating an approximation of the value function.
For iCS, we similarly approximate the nonlinear dynamics about a reference trajectory, but, in contrast, the control updates are found by solving a convex program, which has the benefit of allowing direct consideration of probabilistic constraints at the cost of computation time and restricts the approximation of the dynamics to first order. 

To the best of our knowledge, there are currently no known methods to solve the nonlinear CS problem. Furthermore, in contrast to the existing literature on chance constrained linear CS, we begin our analysis with a continuous stochastic system and describe an exact discretization procedure.

\subsection{Notation and Preliminaries}

For a random vector $z$
on a probability space $(\Omega, \mathscr{F}, \Pr)$,
we denote the expectation of a function $f$ of $z$ as $\E[f(z)]$. The mean value of $z$ is denoted by $\bar{z} := \E(z)$, and the difference from the mean as $\tilde{z} := z - \E(z)$. The covariance of $z$ is denoted by $P_z := \E(\tilde{z} \tilde{z} \t)$.
%and we write $z \sim \mathcal{N}(\bar{z}, P_z)$ to denote that $z$ is normally distributed with mean $\bar{z}$ and covariance $P_z$.
The complement of an event $A \subseteq \Omega$ is denoted $A^c = \Omega \setminus A$, and we use the shorthand $\{\omega \in \Omega : z(\omega) \in B \} = \{ z \in B \}$ to denote events. Dependence of a quantity $y$ on time $t$ is denoted by $y_t$. For a square matrix $A$, we write $A > 0$ $(\geq 0)$ if $A$ is positive (semi-)definite, i.e., $x \t A x > 0$ $(\geq 0)$ for all nonzero real vectors $x$. The set of natural numbers, including zero, is written as $\mathbb{N}_0$, and $\mathbb{N}_+ = \mathbb{N}_0 \setminus \{ 0 \}$. We will denote by $\mathbb{N}_0^m$ the set of natural numbers up to, and including, a positive integer $m$  (similarly for $\mathbb{N}_+^m$).

%%%%%%%%%%%%%%%%%%%%%%%%%%%%%%%%%%%%%%%%%%%%%%%%%%%%%%%%%%%%%%%%%%%%%%%%%%%%%%%%
\section{PROBLEM FORMULATION} \label{sec:problem-formulation}

Consider the nonlinear stochastic differential equation
\begin{equation} \label{eq:original-nonlinear-sde}
	\d x_t = f(x_t, u_t, t) \d t + G_t \d w_t, \quad t \in [t_0, t_f],
\end{equation}
where $x_t \in \Re ^ {n_x}$ is the state, $u_t \in \Re ^ {n_u}$ is the control input, and
$w_t$ is an $n_w$-dimensional standard Brownian motion.
%$\d w_t \in \Re ^{n_w}$ is an increment of standard Brownian motion.
%The functions $f: \Re ^ {n_x} \times \Re ^ {n_u} \times [t_0, t_f] \mapsto \Re ^ {n_x}$ and $G_t : [t_0, t_f] \mapsto \Re ^{n_x \times n_w}$ are deterministic and continuous in time.
At time $t_0$, the state $x_{t_0}$ is assumed to be normally distributed with fixed mean and covariance
\begin{equation} \label{eq:orig-initial-distribution}
	\E(x_{t_0}) = \bar{x}_0, \quad \E( \tilde{x}_{t_0}  \tilde{x}_{t_0} \t ) = P_{x_0}.
\end{equation}
At each time, the state and control are constrained in probability to given convex sets 
\begin{equation} \label{eq:original-chance-constraints}
	\Pr(x_t \in \mathcal{X}_t) \geq 1 - p_{x, t}, \quad \Pr(u_t \in \mathcal{U}_t) \geq 1 - p_{u, t},
\end{equation}
where $0 < p_{x, t} < 1/2$ and $0 < p_{u, t} < 1/2 $ are prescribed maximum probabilities of failure. The constraints (\ref{eq:original-chance-constraints}) are referred to as chance constraints. We wish to find a control that brings the state $x_t$ to a final distribution at time $t_f$
%We wish to find a control $u_t$ that brings the state $x_t$ to a normal distribution at time $t_f$
with given mean and covariance
\begin{equation} \label{eq:orig-terminal-state-constraints}
		\E( x_{t_f} ) = \bar{x}_f, \quad \E( \tilde{x}_{t_f} \tilde{x}_{t_f} \t ) = P_{x_f},
\end{equation}
where $P_{x_f}$ is a given positive-definite matrix, while minimizing the cost functional
\begin{equation} \label{eq:orig-continuous-cost}
	J(u)= \int_{t_0}^{t_f} \big[ \ell(\bar{u}_t, \bar{x}_t) + \E (\tilde{x}_t \t Q_{x,t} \tilde{x}_t + \tilde{u}_t \t Q_{u,t} \tilde{u}_t) \big] \d t.
\end{equation}
Here $Q_{x,t} > 0$ and $Q_{u,t} \geq 0$ are weight matrices, and $\ell$ is an integrable function that is convex in $\bar{u}_t$ and $\bar{x}_t$, and the optimization is performed over the control $u$. In summary, we are interested in solving the following problem.

\begin{problem} \label{problem:orig-nonlinear}
	\textit{Nonlinear Covariance Steering.}
	Find a control $u_t^*$ to minimize the cost (\ref{eq:orig-continuous-cost}) subject to the dynamics (\ref{eq:original-nonlinear-sde}), terminal state constraints (\ref{eq:orig-terminal-state-constraints}), and chance constraints (\ref{eq:original-chance-constraints}).
\end{problem}

In the remainder of this section, we will develop a linear approximation of (\ref{eq:original-nonlinear-sde}) in the neighborhood of a given reference. Then, after discretizing the linearized system, we will focus our analysis on the discrete linear system.

\subsection{Time Normalization and Linearization}

We begin by normalizing the time domain $[t_0, t_f]$ to the unit interval using the dilation coefficient \cite{Szmuk2018}
%The time domain $[t_0, t_f]$ is normalized to the unit interval $[0, 1]$ by the time dilation coefficient
\begin{equation}
	\sigma \defeq t_f - t_0.
\end{equation}
Let $\tau  \defeq (t - t_0) / \sigma \in [0, 1]$ be the normalized time, from which it follows that $\sigma = \d t / \d \tau$. Since the Brownian motion increment $\d w_t$ has variance $\d t$, we scale the diffusion term in (\ref{eq:original-nonlinear-sde}) by $\sqrt{\sigma}$ to obtain $\d w_\tau$ with variance $\d \tau$ (i.e., $\d w_t$ is identically distributed with $\sqrt{\sigma} \d w_\tau$). The time normalized system is then given by
\begin{equation} \label{eq:sde-time-normalized}
	\d x_\tau = \sigma f(x_\tau, u_\tau, \tau) \d \tau + \sqrt{\sigma} G_\tau \d w_\tau, \quad \tau \in [0, 1],
\end{equation}
and the time normalized cost is given by
\begin{equation} \label{eq:continuous-cost-time-normalized}
	J(u)= \sigma \int_{0}^{1} \big[ \ell(\bar{u}_\tau, \bar{x}_\tau) + \E (\tilde{x}_\tau \t Q_{x,\tau} \tilde{x}_\tau + \tilde{u}_\tau \t Q_{u,\tau} \tilde{u}_\tau) \big] \d \tau.
\end{equation}
Next, we linearize (\ref{eq:sde-time-normalized}) about a given reference trajectory $(\hat{x}^i_\tau, \hat{u}^i_\tau)$, where $i \geq 1$ is an index to count successive linearizations. This procedure results in the linear stochastic system
\begin{equation} \label{eq:sde-linearized}
	\d x_\tau \approx (A^i_\tau x_\tau + B^i_\tau u_\tau + r^i_\tau ) \d \tau + \sqrt{\sigma} G_\tau \d w_\tau,
\end{equation}
where
\begin{equation}
	A^i_\tau \defeq \sigma \left. \dfrac{\partial f}{\partial x} \right\vert_{(\hat{x}^i_\tau, \hat{u}^i_\tau)}, \quad
	B^i_\tau \defeq  \sigma \left. \dfrac{\partial f}{\partial u} \right\vert_{(\hat{x}^i_\tau, \hat{u}^i_\tau)}, \quad
\end{equation}
\begin{equation} \label{eq:linearized-d-r}
	r_\tau \defeq \sigma f(\hat{x}^i_\tau, \hat{u}^i_\tau, \tau) - A^i_\tau \hat{x}^i_\tau - B^i_\tau \hat{u}^i_\tau.
\end{equation}

\subsection{Discrete Approximation}

Let $0 = \tau_0 < \tau_1 < \dots < \tau_N = 1$ be a partition of the interval $[0, 1]$, where 
%Partition the interval $[0, 1]$ into $0 = \tau_0 < \tau_1 < \dots < \tau_N = 1$, where
\begin{equation}
	\tau_k \defeq \frac{k}{N}, \quad k \in \N_0^N.
\end{equation}
%It follows that $\d \tau = 1/N$.
Henceforth, we will write $x_k \defeq x_{\tau_k}$ and $u_k \defeq u_{\tau_k}$. We use a zero-order-hold (ZOH) discretization of the control given by
\begin{equation}
	u_\tau = u_k, \quad \tau \in [\tau_k, \tau_{k + 1}), \quad k \in \N_0^{N - 1}.
\end{equation}
Substituting $u_k$ in (\ref{eq:sde-linearized}), we obtain the solution \cite{Fleming1975}
%For simplicity, we will also assume that the noise 
%Furthermore, by assuming a ZOH discretization on the noise, we obtain \textcolor{red}{Fix this.}
%\begin{multline} \label{eq:discrete-full-equation}
%	x_{k + 1} = \Phi(\tau_{k + 1}, \tau_k) x_k + \int_{\tau_k}^{\tau_{k + 1}} \Phi(\tau_{k + 1}, \tau) B^i_\tau \d \tau \, u_k \\
%	+ \int_{\tau_k}^{\tau_{k + 1}} \Phi(\tau_{k + 1}, \tau) d^i_\tau \d \tau \sigma + \int_{\tau_k}^{\tau_{k + 1}} \Phi(\tau_{k + 1}, \tau) r^i_\tau \d \tau \\
%	+ \sqrt{\sigma}  \int_{\tau_k}^{\tau_{k + 1}} \Phi(\tau_{k + 1}, \tau_k) G_\tau \d w_\tau ,
%\end{multline}
\begin{multline} \label{eq:discrete-full-equation}
	x_{k + 1} = \Phi^i(\tau_{k + 1}, \tau_k) x_k \\
	+ \int_{\tau_k}^{\tau_{k + 1}} \Phi^i(\tau_{k + 1}, \tau) ( B^i_\tau \, u_k
	+ r^i_\tau ) \d \tau \\
	+ \sqrt{\sigma}  \int_{\tau_k}^{\tau_{k + 1}} \Phi^i(\tau_{k + 1}, \tau) G_\tau \d w_\tau, \quad k \in \N_0^{N - 1},
\end{multline}
where $\Phi^i(\tau, s)$ is the state transition matrix for system (\ref{eq:sde-linearized}), which satisfies
\begin{equation}
	\dfrac{\partial }{\partial \tau} \Phi^i(\tau, s) = A^i_\tau \Phi^i(\tau, s), \quad \Phi^i(\tau, \tau) = I.
\end{equation}
%Let $G_k$ be the matrix that satisfies the stochastic integral equation
%\begin{equation}
%	\int_{\tau_k}^{\tau_{k + 1}} \Phi(\tau_{k + 1}, \tau) G_\tau \d w_\tau = G_k w_k,
%\end{equation}
%where, for each $k$, $w_k \in \Re ^{n_w}$ are i.i.d normal random vectors with unit variance.
For $k \in \mathbb{N}_0^{N - 1}$, we rewrite (\ref{eq:discrete-full-equation}) as
\begin{equation} \label{eq:discrete-simple-equation}
	x_{k + 1} = A^i_k x_k + B^i_k u_k + r^i_k + \sqrt{\sigma} G^i_k w_k,
\end{equation}
where $w_k \in \Re ^ {n_x}$ are
%i.i.d
independent and identically distributed 
$\mathcal{N}(0, I)$ and where
\begin{subequations}
    \begin{align}
	    A^i_k &\defeq \Phi^i(\tau_{k + 1}, \tau_k), \\
	    B^i_k &\defeq \int_{\tau_k}^{\tau_{k + 1}} \Phi^i(\tau_{k + 1}, \tau) B^i_\tau \d \tau, \\
		r^i_k &\defeq \int_{\tau_k}^{\tau_{k + 1}} \Phi^i(\tau_{k + 1}, \tau) r^i_\tau \d \tau.
%  	 	G^i_k &\defeq  \left( \int_{\tau_k}^{\tau_{k + 1}} \Phi^i(\tau_{k + 1}, \tau) G_\tau G_\tau \t \Phi ^{i \addt} (\tau_{k + 1}, \tau) \d \tau \right)^{\nicefrac{1}{2}} \negthickspace\negthickspace\negthickspace.\negthickspace
    \end{align}
\end{subequations}
The stochastic integral in (\ref{eq:discrete-full-equation}) is a zero-mean Gaussian random vector with covariance
\begin{equation}
	\Sigma = \int_{\tau_k}^{\tau_{k + 1}} \Phi^i(\tau_{k + 1}, \tau) G_\tau G_\tau \t \Phi ^{i \addt} (\tau_{k + 1}, \tau) \d \tau ,
\end{equation}
and therefore the matrix $G^i_k$ is selected such that $G^i_k w_k \in \Re ^{n_x}$ is distributed $\mathcal{N}(0, \Sigma)$.  Since $\Sigma$ may have rank greater than $n_w$ (e.g., $\Sigma$ is full rank if $(A^i_\tau, G_\tau)$ is controllable), the matrix $G^i_k$ must be $n_x \times n_x$.
Furthermore, since the solution to $G^i_k G^{i \addt}_k = \Sigma$ is not unique, there exist a continuum of coefficient matrices $G^i_k$ such that the resulting state processes are equidistributed.

%The matrix $G^i_k \in \Re^{n_x \times n_x}$ is selected so that the vector $G^i_k w_k$ is $\mathcal{N}(0, \Sigma)$, where
%\begin{eqnarray}
%	\Sigma = \int_{\tau_k}^{\tau_{k + 1}} \Phi^i(\tau_{k + 1}, \tau) G_\tau G_\tau \t \Phi ^{i \addt} (\tau_{k + 1}, \tau) \d \tau .
%\end{eqnarray}
%It follows that $G^i_k$ is not unique, but any $G^i_k$ such that $G^i_k G^{i \addt}_k = \Sigma$ will result in an distributed process $x_k$. Furthermore, note that
%It follows that $G^i_k$ is not unique, but we (ref).
%that is, $G^i_k G^{i \addt}_k = \Sigma$, where $\Sigma$ is an $n_x$ by $n_x$ positive semidefinite matrix. It follows that $G^i_k$ is unique modulo unitary transformation. TODO
%\begin{remark}
The discrete formulation (\ref{eq:discrete-full-equation}) is an exact representation of the linear system (\ref{eq:sde-linearized}) with ZOH control. However, since the previous integrals may be difficult to compute, the following first-order approximation is commonly used:
\begin{flalign}
    \begin{aligned}
    &A^i_k \defeq A^i_{\tau_k} \d \tau + I, &\quad B^i_k &\defeq B^i_{\tau_k} \d \tau, \\
    &G^i_k \defeq \sqrt{\d \tau} G_{\tau_k}, &\quad r^i_k &\defeq r^i_{\tau_k} \d \tau.
    \end{aligned}
\end{flalign}

%\end{remark}
We remark that the method of discretization chosen does not affect any of the following discussion.

Following \cite{Goldshtein2017, Okamoto2018b, Okamoto2018}, we will rewrite the discrete system (\ref{eq:discrete-simple-equation}) as a single linear equation.
%
%To this end, we define the transition matrices from step $k_0$ to step $k_1$ as
%$A^i_{k_1, k_0} := A^i_{k_1} A^i_{k_1- 1} \cdots A^i_{k_0}$,
%$B^i_{k_1, k_0} := A^i_{k_1, k_0 + 1} B_{k_0}$, 
%$G_{k_1, k_0} := A^i_{k_1, k_0 + 1} G_{k_0}$,
%and similarly, we define the vectors
%$d^i_{k_1, k_0} := A^i_{k_1, k_0 + 1} d_{k_0}$,
%$r^i_{k_1, k_0} := A^i_{k_1, k_0 + 1} r_{k_0}$.
Define the state transition matrix from step $k_0$ to step $k_1$ as
\begin{equation}
	A^i_{k_1, k_0} \defeq \begin{cases}
		A^i_{k_1} A^i_{k_1- 1} \cdots A^i_{k_0}, & k_1 \geq k_0, \\
		I, & k_1 < k_0,
	\end{cases}
\end{equation}
%\begin{equation}
%	A^i_{k_1, k_0} \defeq A^i_{k_1} A^i_{k_1- 1} \cdots A^i_{k_0},
%\end{equation}
and the corresponding transitions for the control, noise, and affine terms as
\begin{equation}
	B^i_{k_1, k_0} \defeq A^i_{k_1, k_0 + 1} B_{k_0},  \quad G^i_{k_1, k_0} \defeq A^i_{k_1, k_0 + 1} G^i_{k_0},
\end{equation}
\begin{equation}
	r^i_{k_1, k_0} \defeq A^i_{k_1, k_0 + 1} r_{k_0}.
\end{equation}
Define concatenated control and disturbance vectors at step $k$ as
\begin{align}
	U_k &\defeq \begin{bmatrix}
		u_0 \t & u_1 \t & \cdots & u_k \t
	\end{bmatrix} \t \in \Re ^{(k + 1)n_u}, \\
	W_k &\defeq \begin{bmatrix}
		w_0 \t & w_1 \t & \cdots & w_k \t
	\end{bmatrix} \t \in \Re ^{(k + 1)n_x}.
\end{align}
Then the state at step $k$ can be written as
\begin{equation} \label{eq:discrete-step-k}
	x_k = \bar{A}^i_k x_0 + \bar{B}^i_k U_{k - 1} + \bar{r}^i_k 1_k + \sqrt{\sigma} \bar{G}^i_k W_{k - 1},
\end{equation}
where $1_k \in \Re^k$ is a column vector of ones, $\bar{A}^i_k \defeq A^i_{k - 1, 0}$, and
\begin{align}
	\bar{B}^i_k &\defeq \begin{bmatrix}
		B^i_{k-1,0} & B^i_{k-1,1} & \cdots & B^i_{k-1,k-2} & B^i_{k - 1}
	\end{bmatrix}, \\
	\bar{G}^i_k &\defeq \begin{bmatrix}
			G^i_{k-1,0} & G^i_{k-1,1} & \cdots & G^i_{k-1,k-2} & G^i_{k - 1}
		\end{bmatrix}, \\
	\bar{r}^i_k &\defeq \begin{bmatrix}
		r^i_{k-1,0} & r^i_{k-1,1} & \cdots & r^i_{k-1,k-2} & r^i_{k - 1}
	\end{bmatrix}.
\end{align}
In terms of the concatenated state vector $X \defeq \begin{bmatrix}
	x_0 \t & \cdots & x_N \t
\end{bmatrix} \t \in \Re ^{(N + 1)n_x}$, control vector $U \defeq U_{N - 1} \in \Re ^{N n_u}$, and disturbance vector $W \defeq W_{N - 1} \in \Re ^{N n_x}$, the system dynamics are written as the matrix equation
\begin{equation} \label{eq:main-linear-vector-equation}
	X = \Amat^i x_0 + \Bmat^i U + R^i + \sqrt{\sigma} \Gmat^i W.
\end{equation}
The matrices $\Amat^i, \Bmat^i, R^i$, and $\Gmat^i$ are formed by appropriately stacking the terms from (\ref{eq:discrete-step-k}).

%Before proceeding, we briefly review our progress up to this point. Our objective is to control the nonlinear stochastic system (\ref{eq:original-nonlinear-sde}) defined on the time interval $[t_0, t_f]$, where the final time $t_f$ is free. We normalized the time interval to $[0, 1]$, and then linearized about a reference time dilation and trajectory to obtain (\ref{eq:sde-linearized}). The resulting linear system was discretized and rewritten as the single linear equation (\ref{eq:main-linear-vector-equation}), which will serve as the foundation of our solution approach in the following sections.

Let $\mathcal{Q}_x$ and $\mathcal{Q}_u$ be block-diagonal state and control cost weight matrices with entries corresponding to the continuous weights $Q_{x,t}$ and $Q_{u,t}$ from (\ref{eq:orig-continuous-cost}):
\begin{equation}
	\mathcal{Q}_x \defeq \mathrm{blkdiag}(Q_{x, \tau_0}, \dots, Q_{x, \tau_{N - 1}}, 0_{n_x}),
\end{equation}
\begin{equation}
	\mathcal{Q}_u \defeq \mathrm{blkdiag}(Q_{u, \tau_0}, \dots, Q_{u, \tau_{N - 1}}).
\end{equation}
The quadratic state cost at step $N$ is neglected, since the terminal state is fixed. 
Letting $E_k = \begin{bmatrix}
	0_{n_x, k n_x}, I_{n_x}, 0_{n_x, (N - k) n_x}
\end{bmatrix}$
and $E^u_k = \begin{bmatrix}
	0_{n_u, k n_u}, I_{n_u}, 0_{n_u, (N - k - 1) n_u}
\end{bmatrix}$, such that $x_k = E_k X$ and $u_k = E^u_k U$,
%and $R_k = \begin{bmatrix}
%	0_{n_u, k n_u}, I_{n_u}, 0_{n_u, (N - k - 1) n_u}
%\end{bmatrix}$.
the continuous time cost functional (\ref{eq:continuous-cost-time-normalized}) can be rewritten in terms of the linearized system as
\begin{equation} \label{eq:cost-discrete}
	J(U) \approx \frac{\sigma}{N} \bigg[ \sum_{k = 0}^{N - 1} \ell(E^u_k \bar{U}, E_k \bar{X}) + \E( \tilde{X} \t \mathcal{Q}_x \tilde{X} + \tilde{U} \t \mathcal{Q}_u \tilde{U} \t ) \bigg],
\end{equation}
and the boundary conditions (\ref{eq:orig-initial-distribution}) and (\ref{eq:orig-terminal-state-constraints}) as
\begin{subequations} \label{eq:discrete-boundary-conditions}
	\begin{alignat}{3}
		E_0 \bar{X} &= \bar{x}_0, \quad &E_0 \E (\tilde{X} \tilde{X} \t) E_N \t &= P_{x_0}, \\
		\label{eq:discrete-boundary-conditions-final}
		E_N \bar{X} &= \bar{x}_f, \quad &E_N \E (\tilde{X} \tilde{X} \t) E_N \t &= P_{x_f}.
	\end{alignat}
\end{subequations}
The chance constraints (\ref{eq:original-chance-constraints}), enforced at each time step $k$, are written as
\begin{subequations} \label{eq:discrete-chance-constraints}
	\begin{align}
		\Pr(E_k X \in \mathcal{X}_{\tau_k}) &\geq 1 - p_{x, \tau_k}, \\
		\Pr(E^u_k U \in \mathcal{U}_{\tau_k}) &\geq 1 - p_{u, \tau_k}.
	\end{align}
\end{subequations}

In summary, we have approximated the continuous time, nonlinear stochastic system (\ref{eq:original-nonlinear-sde}) by the discrete, linear stochastic system (\ref{eq:main-linear-vector-equation}). Problem \ref{problem:orig-nonlinear} can be accordingly restated in terms of this approximate system as follows.

\begin{problem}
	Find the control sequence $U^*$ that minimizes (\ref{eq:cost-discrete}) subject to the dynamics (\ref{eq:main-linear-vector-equation}), boundary conditions (\ref{eq:discrete-boundary-conditions}), and chance constraints (\ref{eq:discrete-chance-constraints}).
\end{problem}

\begin{remark}
	In the discrete-time formulation, the chance constraints are only enforced at the discrete times $\tau_k$, and therefore the original constraints (\ref{eq:original-chance-constraints}) may be violated for some $\tau \in (\tau_k, \tau_{k + 1})$. Constraint violation in this interval is likely when the discretization is too coarse.
\end{remark}

\section{COVARIANCE STEERING} \label{sec:covariance-steering}

For the remainder of this paper, we will restrict the control law to be of the form \cite{Okamoto2018}
\begin{equation}
	u_k = v_k + K_k y_k,
\end{equation}
were $v_k \in \Re ^ {n_u}$ is a feedforward control, $K_k \in \Re ^{n_u \times n_x}$ is a feedback gain matrix, and $y_k \in \Re ^ {n_x}$ is a zero-mean random process given by
\begin{equation}
	y_{k + 1} = A^i_k y_k + \sqrt{\sigma} G^i_k w_k, \quad y_0 = x_0 - \bar{x}_0.
\end{equation}
In vector notation, we have
\begin{equation}
	Y = \Amat^i y_0 + \sqrt{\sigma} \Gmat^i W,
\end{equation}
and thus,
\begin{equation}
	U = V + KY = V + K(\Amat^i y_0 + \sqrt{\sigma} \Gmat^i W),
\end{equation}
where the block feedback matrix $K  \in \Re ^{N n_u \times (N + 1) n_x}$ is given by
\begin{equation}
	K \defeq \begin{bmatrix}
		\mathrm{blkdiag}(K_0, \dots, K_{N - 1}) & 0_{N n_u, n_x}
	\end{bmatrix}.
\end{equation}
Substituting the control into the state equation, we obtain the expressions for the mean and deviation states as
\begin{align}
	\label{eq:Xbar}
	\bar{X} &\defeq \E (X) = \Amati \bar{x}_0 + \Bmati V + R^i, \\
	\tilde{X} &\defeq X - \E(X) \nonumber \\
	&\,= \Amati y_0 + \Bmati K (\Amati y_0 + \sqrt{\sigma} \Gmat^i W) + \sqrt{\sigma} \Gmat^i W \nonumber \\
	&\,= (I + \Bmati K)(\Amati y_0 + \sqrt{\sigma} \Gmat^i W) .
\end{align}
Similarly for the control, we obtain
\begin{align}
	\bar{U} &\defeq \E (U) = V, \\
	\tilde{U} &\defeq U - \E(U) = K (\Amati y_0 + \sqrt{\sigma} \Gmat^i W).
\end{align}
It follows that the state and control covariances, in terms of the covariance of the process $y_k$, are given as
%\begin{align}
%	\label{eq:matrix-state-cov}
%	\mathcal{P}_x &\defeq \E (\tilde{X} \tilde{X} \t) \\
%	&\,= (I + \Bmati K)(\Amat^i P_{x_0} \Amatit + \sigma \Gmat \Gmat \t)(I + \Bmati K) \t, \\
%	\label{eq:matrix-control-cov}
%	\mathcal{P}_u &\defeq \E (\tilde{U} \tilde{U} \t ) = K (\Amati P_{x_0} \Amatit + \sigma \Gmat \Gmat \t) K \t .
%\end{align}
\begin{alignat}{3}
	\mathcal{P}_y &\defeq \E(Y Y\t) &&= \Amat^i P_{x_0} \Amatit + \sigma \Gmat ^i \Gmat ^{i \addt}, \\
 	\label{eq:matrix-state-cov}
	\mathcal{P}_x &\defeq \E (\tilde{X} \tilde{X} \t) &&= (I + \Bmati K) \mathcal{P}_y (I + \Bmati K) \t, \\
	\label{eq:matrix-control-cov}
	\mathcal{P}_u &\defeq \E (\tilde{U} \tilde{U} \t ) &&= K \mathcal{P}_y K \t .
\end{alignat}
%Let $L$ be the mean cost function of the concatenated mean control $\bar{U}$ and mean state $\bar{X}$ corresponding to the running cost function $\ell$ in (\ref{eq:orig-continuous-cost}):
%\begin{equation}
%	L(\bar{U}) = \sum_{k = 0}^{N - 1} \ell(E^u_k \bar{U}, E_k \bar{X}),
%\end{equation}
Substituting (\ref{eq:matrix-state-cov}) and (\ref{eq:matrix-control-cov}) into the cost function (\ref{eq:cost-discrete}) and simplifying, we obtain
%\begin{multline} \label{eq:expanded-matrix-cost-funciton}
%	J(V, K, \sigma) = L(V, \bar{X}) \\ + \tr \big\{ \big[ (I + \Bmat^i K) \t \mathcal{Q}_x (I + \Bmat^i K) + K \t \mathcal{Q}_u K \big] \times \\\big(\Amat^i P_{x_0} \Amatit + \sigma \Gmat \Gmat \t \big) \big\}.
%\end{multline}
\begin{multline} \label{eq:expanded-matrix-cost-funciton}
	J(V, K) = \frac{\sigma}{N} \bigg[ L(V) + \tr \big\{ \big[ (I + \Bmat^i K) \t \mathcal{Q}_x (I + \Bmat^i K) \\ + K \t \mathcal{Q}_u K \big] \mathcal{P}_y  \big\} \bigg],
\end{multline}
%\begin{multline} \label{eq:expanded-matrix-cost-funciton}
%	J(V, K) = \sigma \big[L(V) + S(K) \big],
%\end{multline}
where
\begin{equation}
	L(V) \defeq \sum_{k = 0}^{N - 1} \ell \big(E^u_k V, E_k ( \Amati \bar{x}_0 + \Bmati V + R^i) \big).
\end{equation}
%\begin{equation}
%	S(K) \defeq \tr \big\{ \big[ (I + \Bmat^i K) \t \mathcal{Q}_x (I + \Bmat^i K) + K \t \mathcal{Q}_u K \big] \mathcal{P}_y  \big\}
%\end{equation}
\subsection{Endpoint Constraints}

Substituting (\ref{eq:Xbar}) into (\ref{eq:discrete-boundary-conditions-final}), we obtain the equality constraint on the final mean state
%In terms of $E_N = \begin{bmatrix}
%	0 & 0 & \cdots & I
%\end{bmatrix} \in \Re ^{n_x \times (N + 1) n_x}$, the terminal mean state constraint is given by the equality constraint
\begin{equation} \label{eq:final-mean-equality}
	h(V) := E_N \big( \Amat^i \bar{x}_0 + \Bmat^i V + R^i \big) - \bar{x}_f = 0.
\end{equation}
Since the equality constraint $E_N \mathcal{P}_x E_N \t = P_{x_f}$ is not convex in $K$, and since in practice a smaller than anticipated state covariance is acceptable, we instead enforce the relaxed inequality constraint \cite{Bakolas2016}
%\begin{equation}
%	E_N (I + \Bmati K)(\Amati P_{x_0} \Amatit + \sigma \Gmat \Gmat \t)(I + \Bmati K) \t E_N \t \leq P_{x_f},
%\end{equation}
\begin{equation}
	E_N (I + \Bmati K) \mathcal{P}_y (I + \Bmati K) \t E_N \t \leq P_{x_f},
\end{equation}
which is convex in $K$. This constraint may be equivalently stated in the more standard form \cite{Okamoto2018b}
%\begin{multline} \label{eq:final-covariance-inequality}
%	g(K, \sigma) := \Vert (\Amati P_{x_0} \Amatit + \sigma \Gmat \Gmat \t)^{1/2}(I + \Bmati K) \t \times \\ E_N \t P_{x_f} ^{-1/2} \Vert_2 - 1 \leq 0.
%\end{multline}
\begin{multline} \label{eq:final-covariance-inequality}
	g(K) := \Vert \mathcal{P}_y^{1/2}(I + \Bmati K) \t E_N \t P_{x_f} ^{-1/2} \Vert_2 - 1 \leq 0.
\end{multline}
\subsection{Chance Constraints}

Assume that at each time step the convex regions $\mathcal{X}_k \defeq \mathcal{X}_{\tau_k}$ and $\mathcal{U}_k \defeq \mathcal{U}_{\tau_k}$ can be represented by the finite intersection of half spaces
%\begin{align}
%	\mathcal{X}_k &= \bigcap_{m = 1}^{M_x} \{ x \in \Re^{n_x} : a_{\tau, m} \t x \leq \alpha_{\tau, m} \}, \\
%	\mathcal{U}_k &= \bigcap_{m = 1}^{M_u} \{ u \in \Re ^{n_u} : b_{\tau, m} \t u \leq \beta_{\tau, m} \},
%\end{align}
\begin{align}
	\mathcal{X}_k = \bigcap_{m = 1}^{M_x} \mathcal{X}_{k, m}, \quad \mathcal{U}_k = \bigcap_{m = 1}^{M_u} \mathcal{U}_{k, m},
\end{align}
where $\mathcal{X}_{k, m} \defeq \{ x \in \Re^{n_x} : \bar{a}_{k, m} \t x \leq \alpha_{k, m} \}$ and $\mathcal{U}_{k, m} \defeq \{ u \in \Re ^{n_u} : \bar{b}_{k, m} \t u \leq \beta_{k, m} \}$ are given in terms of the vectors $\bar{a}_{k, m} \in \Re ^{n_x}$, $\bar{b}_{k, m} \in \Re ^{n_u}$ and scalars $\alpha_{k, m}, \beta_{k, m} \in \Re$.
%
%As part of our discretization scheme, we concatenated the states and controls into vectors in $\Re^{(N + 1) n_x}$ and $\Re^{N n_u}$. Therefore, in order to enforce these constraints at each step $k$, we will rewrite the constraint sets $\mathcal{X}_\tau \subset \Re ^{n_x}$ and $\mathcal{U}_\tau \subset \Re ^{n_u}$ as single sets $\mathcal{X} \subset \Re ^{(N + 1) n_x}$ and $\mathcal{U} \subset \Re ^{N n_u}$ with the property
%\begin{align}
%	X \in \mathcal{X} &\implies x_k = E_k X \in \mathcal{X}_{\tau_k}, \\
%	U \in \mathcal{U} &\implies u_k = E^u_k X \in \mathcal{U}_{\tau_k}.
%\end{align}
%
By subadditivity of probability, we have
\begin{align}
	\Pr (x_k \in \mathcal{X}_k^c) &= \Pr \bigg(x_k \in \bigcup_{m = 1}^{M_x} \mathcal{X}^c_{k, m} \bigg) \leq \sum_{m = 1}^{M_x} \Pr(x_k \in \mathcal{X}^c_{k, m} ). \label{eq:subadditivity-state-ineq}
\end{align}
It follows that if $\Pr(x_k \in \mathcal{X}^c_{k, m} ) \leq p^x_{k, m}$ for a set of positive numbers $\{ p^x_{k, m} \}$ that sum over the index $m$ to less than $p_{x, k}$, then $\Pr (x_k \in \mathcal{X}_k^c) \leq p_{x, k}$ \cite{Nemirovski2006, Blackmore2009}.
In terms of the concatenated state and control vectors, and since $x_k = E_k X$ and $u_k = E^u_k U$, 
the events
$\{ x_k \in \mathcal{X}_k \} \subset \Re^{n_x}$ and $\{ E_k X \in \mathcal{X}_k \} \subset \Re^{(N + 1)n_x}$ have the same probability.
%it follows from the definition of a product measure that
%\begin{equation}
%	\Pr(x_k \in \mathcal{X}_k) = \Pr(E_k X \in \mathcal{X}_k)
%\end{equation}
Therefore, when relabeling indices of the inequality constraints according to
%Relabeling the indices of the inequalities according to
%For the discrete system, we will consider these sets at each time step $k$, where
\begin{alignat}{2}
	a_m \t E_k &= \bar{a}_{k, m} \t , \quad && m \in \mathbb{N}_+^{M_x}, \; k \in \mathbb{N}_0^{N}, \\
	b_m \t E^u_k &= \bar{b}_{k, m} \t , \quad && m \in \mathbb{N}_+^{M_u}, \; k \in \mathbb{N}_0^{N - 1},
\end{alignat}
if $\{ p^x_{m, k} \}$ and $\{ p^u_{m, k} \}$ are given sets of positive numbers that satisfy, for each $k$,
\begin{equation} \label{eq:probability-weight-sum}
	\sum_{m = 1}^{M_x} p^x_{m, k} \leq p_{x, k}, \quad \sum_{m = 1}^{M_u} p^u_{m, k} \leq p_{u, k},
\end{equation}
then, from (\ref{eq:subadditivity-state-ineq}),
%\begin{multline} \label{eq:probability-implication}
%	\Pr(a_m \t E_k X \leq \alpha_{k, m}) \geq 1 - p^x_{m, k} \quad \forall m \in \mathbb{N}_+^{M_x} \\
%		\implies \Pr (x_k \in \mathcal{X}_k) \geq 1 - p_{x, k}.
%\end{multline}
\begin{equation}
	\Pr(a_m \t E_k X \leq \alpha_{k, m}) \geq 1 - p^x_{m, k}, \quad m \in \mathbb{N}_+^{M_x},
\end{equation}
it follows that
\begin{equation}
	\Pr (x_k \in \mathcal{X}_k) \geq 1 - p_{x, k}.
\end{equation}
The same construction applies to the control sets $\mathcal{U}_k$.

Next, we formulate the chance constraint $\Pr(a_m \t E_k X \leq \alpha_{k, m})$ into a deterministic expression of the control variables. From (\ref{eq:main-linear-vector-equation}) it follows that $X$ is normally distributed and hence $a_m \t E_k X$ is a scalar normal random variable with mean $a_m \t E_k \bar{X}$ and covariance $a_m \t E_k \mathcal{P}_x E_k \t a_m$. It follows that
\begin{equation}
	\Pr(a_m \t E_k X \leq \alpha_{k, m}) = \ncdf \bigg( \frac{\alpha_{k, m} - a_m \t E_k \bar{X}}{\sqrt{a_m \t E_k \mathcal{P}_x E_k \t a_m}} \bigg),
\end{equation}
where $\ncdf$ is the cumulative normal distribution function. 
Therefore, the chance constraint $\Pr(a_m \t E_k X \leq \alpha_{k, m}) \geq 1 - p_{x, k}$ can be equivalently written as
\begin{equation}
	a_m \t E_k \bar{X} - \alpha_{k, m} + \ncdf \inv (1 - p_{x, k}) \big\Vert \big( \mathcal{P}_x E_k \t a_m \big)^{1/2} \big\Vert \leq 0.
\end{equation}
Putting it all together, if (\ref{eq:probability-weight-sum}) holds, and if
\begin{subequations} \label{eq:chance-constraint-final}
	\begin{multline}
		\label{eq:chance-constraint-final-state}
		c^x_{m, k}(K, V) \defeq a_m \t E_k (\Amat^i \bar{x}_0 + \Bmat^i V + R^i)- \alpha_m  \\ + \mathrm{cdfn} \inv (1 - p^x_{m,k}) \big\Vert \mathcal{P}_y^{1/2} (I + \Bmat^ i K) \t E_k \t a_m \big\Vert \leq 0,
	\end{multline}
	and
	\begin{multline}
		\label{eq:chance-constraint-final-control}
		c^u_{m, k}(K, V) \defeq b_{m} \t V - \beta_m \\ + \mathrm{cdfn} \inv (1 - p^u_{m,k}) \big\Vert \mathcal{P}_y^{1/2} K \t E_k^{u \addt} b_m \big\Vert \leq 0,
	\end{multline}
\end{subequations}
then we ensure that the chance constraints (\ref{eq:discrete-chance-constraints}) will be satisfied.
%\begin{remark}
%	Fixing the  constants $p^x_{m, k}$ and $p^u_{m, k}$ may lead to an overly conservative solution. However, optimizing over these values, subject to the constraint (\ref{eq:probability-weight-sum}), causes the chance constraints to be non-convex \cite{Nemirovski2006}. For simplicity, we use the heuristic
%	\begin{equation}
%		p^x_{m, k} = p_{x, k} / M_x, \quad p^u_{m, k} = p_{u, k} / M_u.
%	\end{equation}
%	The problem of optimizing these values subject to the original chance constraints is referred to as the risk allocation problem \cite{Ono2008, Blackmore2009, Ono2008b}. 
%\end{remark}
\begin{remark}
	This work assumes that $\{ p^x_{m, k} \}$ and $\{ p^u_{m, k} \}$ are given sets of positive numbers that satisfy~(\ref{eq:probability-weight-sum}). 
	This assumption allows~(\ref{eq:chance-constraint-final}) to be convex. 
	Otherwise,~(\ref{eq:chance-constraint-final}) becomes non-convex, and we need to consider an optimal risk allocation problem. 
	Several approaches have been proposed, such as~\cite{Ono2008b,vitus2011closed}, to handle the risk allocation problem.
	In addition, the authors of~\cite{ma2012fast} used a primal-dual interior point method to find an optimal risk allocation.
\end{remark}
We are now ready to restate the covariance steering problem as a deterministic, finite dimensional optimization problem.

\begin{problem} \label{problem:linear-nonconvex}
	\textit{Linear Covariance Steering.} Find $K^*$ and $V^*$ that minimize the cost (\ref{eq:expanded-matrix-cost-funciton}) subject to the terminal state constraints (\ref{eq:final-mean-equality}) and (\ref{eq:final-covariance-inequality}) and the chance constraints (\ref{eq:chance-constraint-final}).
\end{problem}

%\begin{remark}
%	Problem \ref{problem:linear-nonconvex} is convex for a fixed $\sigma$, however, when optimizing over $\sigma$, the cost (\ref{eq:expanded-matrix-cost-funciton}), the terminal state covariance constraint (\ref{eq:final-covariance-inequality}), and the chance constraints (\ref{eq:chance-constraint-final}) are not convex.
%\end{remark}

% \section{SUCCESSIVE CONVEXIFICATION}
\section{ITERATIVE COVARIANCE STEERING}
\label{sec:successive-convexification}

In the previous sections, we have locally approximated the continuous time, nonlinear system (\ref{eq:original-nonlinear-sde}) with the discrete linear system (\ref{eq:main-linear-vector-equation}), and we have restated the cost function and constraints in terms of the discrete linear system as functions of a feedfoward control $V$ and feedback gain $K$.
%Next, we further approximate the linear discrete system so that the optimal controls may be found by convex programming.
We will search for solutions to the original nonlinear system by successively solving this approximate convex problem, where the optimal controls from each successive problem are used to propagate trajectories of the \textit{nonlinear} system to obtain references for the next linearization step. This method is referred to in the literature as successive convexification \cite{Szmuk2016, Mao2016}.

\subsection{Stochastic Trust Region}

The linear approximation of the system dynamics is only valid in a neighborhood around the reference trajectory, so care must be taken to ensure that the optimal controls for the linear problem are relevant to the nonlinear problem. For this reason, variations in the state and the control from the previous solution are bounded inside a \textit{trust region} \cite{Szmuk2016}.
In this paper, we are successively approximating a stochastic system, and since the state of a system with Brownian noise is unbounded, we must define a stochastic trust region instead as follows
\begin{subequations} \label{eq:trust-region-constraint}
	\begin{equation} % \label{eq:trust-region-constraint}
		\Pr \big( \norm{\hat{x}^i_k - x_k}_1 \leq \Delta_x^i \big) \geq 1 - p^x_\mathrm{tr}, \quad k \in \N_0^{N},
	\end{equation}
	\begin{equation} % \label{eq:trust-region-constraint}
		\Pr \big( \norm{\hat{u}^i_k - u_k}_1 \leq \Delta_u^i \big) \geq 1 - p^u_\mathrm{tr}, \quad k \in \N_0^{N - 1},
	\end{equation}
\end{subequations}
where $p^x_\mathrm{tr}$, $p^u_\mathrm{tr}$, $\Delta_x^i$, and $\Delta_u^i$ are user-defined limits. Constraints in the 1-norm can be represented by $2 n_x$ or $2 n_u$ inequality constraints for the state or control, respectively, using (\ref{eq:chance-constraint-final}).
%For this reason, the mean controls $V$ and state $\bar{X}$ are constrained to balls around the reference trajectory
%\begin{equation} \label{eq:trust-region-constraint}
%	\norm{\hat{x}^i_k - E_k \bar{X}} \leq \Delta_x^i, \quad \norm{\hat{u}^i_k - E_k^u V} \leq \Delta_u^i,
%\end{equation}
As a consequence of these trust region constraints, if the reference trajectory is sufficiently far away from the terminal constraint, then the problem may become infeasible. For these situations, which are most likely encountered when initializing the problem, we relax the hard constraint (\ref{eq:final-mean-equality}) on the terminal state mean to the soft constraint
\begin{equation} \label{eq:terminal-mean-relaxed}
	\norm{ E_N (\Amat^i \bar{x}_0 + \Bmat^i V + R^i) - \bar{x}_f} \leq \eta_{x_f},
\end{equation}
with a corresponding term $\eta_{x_f} w_{x_f}$ added to the cost, where $w_{x_f}$ is a user-defined weight. This constraint may be replaced with the hard constraint (\ref{eq:final-mean-equality}) when the reference trajectory $\hat{x}^i$ is sufficiently close to the terminal constraint. 
In the case (\ref{eq:terminal-mean-relaxed}) is active, we use the augmented cost function given by
\begin{equation} \label{eq:cost-convex-approx}
	\mathcal{J}(V, K, \eta_{x_f}) = J(V, K) + \eta_{x_f} w_{x_f}.
\end{equation}
In summary, we have modified Problem \ref{problem:linear-nonconvex} to the following convex optimization problem.

\begin{problem} \label{prob:convex-cs-subproblem}
	\textit{iCS Convex Subproblem.} Find $K^*$ and $V^*$ that minimize the cost (\ref{eq:cost-convex-approx}) subject to the terminal state constraints (\ref{eq:chance-constraint-final}) and (\ref{eq:final-mean-equality}) (or (\ref{eq:terminal-mean-relaxed}) if the reference trajectory is sufficiently far from the target), the chance constraints (\ref{eq:chance-constraint-final}), and the trust region constraints (\ref{eq:trust-region-constraint}).
\end{problem}

This problem is solved successively in order to find a solution to Problem \ref{problem:orig-nonlinear} using the iCS algorithm presented in Procedure \ref{algo:iCS}.

\begin{algorithm}
% 	\caption{Successive Convexificiation}\label{algo:active}
	\caption{Iterative Covariance Steering (iCS)}\label{algo:iCS}
	\begin{algorithmic}[1]
		\Require Initial guess $\hat{u}_k^1, \hat{K}_k^1$
		\Ensure Optimal control $\bar{u}_k^*$ and $K_k^*$
%		\Procedure{successive\_convexification}{$\hat{\sigma}^1, \hat{u}^1_\tau$}
			\For {$i = 1$ to $i_\mathrm{max}$}
%				\State Propagate the mean of the nonlinear system (ref) with control $\hat{\sigma}^i$ and time dilation $\hat{\sigma}^i$ to obtain $\bar{x}$
				\State Propagate nonlinear mean dynamics with $\hat{u}_k^i, \hat{K}_k^i$
%				\State Forward simulate the nonlinear system with control $\hat{\sigma}^i$ and time dilation $\hat{\sigma}^i$ to obtain $\bar{x}$
				\State $\hat{x}_k^i \leftarrow \bar{x}_k$
				\State Linearize about $(\hat{x}^i, \hat{u}^i)$% using (\ref{eq:sde-linearized})
				\State Discretize% the linearized system
				\State Solve problem (\ref{prob:convex-cs-subproblem}) to obtain $V^*, K^*$
				\State Reshape $\bar{u}^*_k \leftarrow V^*$, $K_k^* \leftarrow K^*$
				\If{$\max_{k \in \N_0^{N - 1}} \norm{\bar{u}^*_k - \hat{u}^i_k} \leq \mathrm{tol}$}
					\State \Return $\bar{u}_k^*, K_k^*$
				\Else
					\State $\hat{u}_k^{i + 1} \leftarrow \bar{u}_k^*$, $\hat{K}_k^{i + 1} \leftarrow K_k^*$
				\EndIf
			\EndFor
		\State \Return Convergence not met
%			\State \Return $V^*, K^*, \sigma^*$
%		\EndProcedure
	\end{algorithmic}
\end{algorithm}

\begin{remark}
    The nonlinear mean dynamics can be propagated through Monte Carlo, which can be parallelized. 
    In the case when computational resources are limited, we can approximate $\E[f(x_t, u_t, t)] \approx f(\bar{x}_t, \bar{u}_t, t)$ so that the mean state evolves according to
    \begin{equation}
        \Dot{\bar{x}}_t = f(\bar{x}_t, \bar{u}_t, t).
    \end{equation}
    In this case, the mean state can be estimated by integrating a single trajectory.
\end{remark}

\section{NUMERICAL EXAMPLE} \label{sec:numerical-example}

\begin{figure}
	\centering
	% trim={<left> <lower> <right> <upper>}
	\includegraphics[trim={0.15in 0.05in 0.1in 0.08in}, clip]{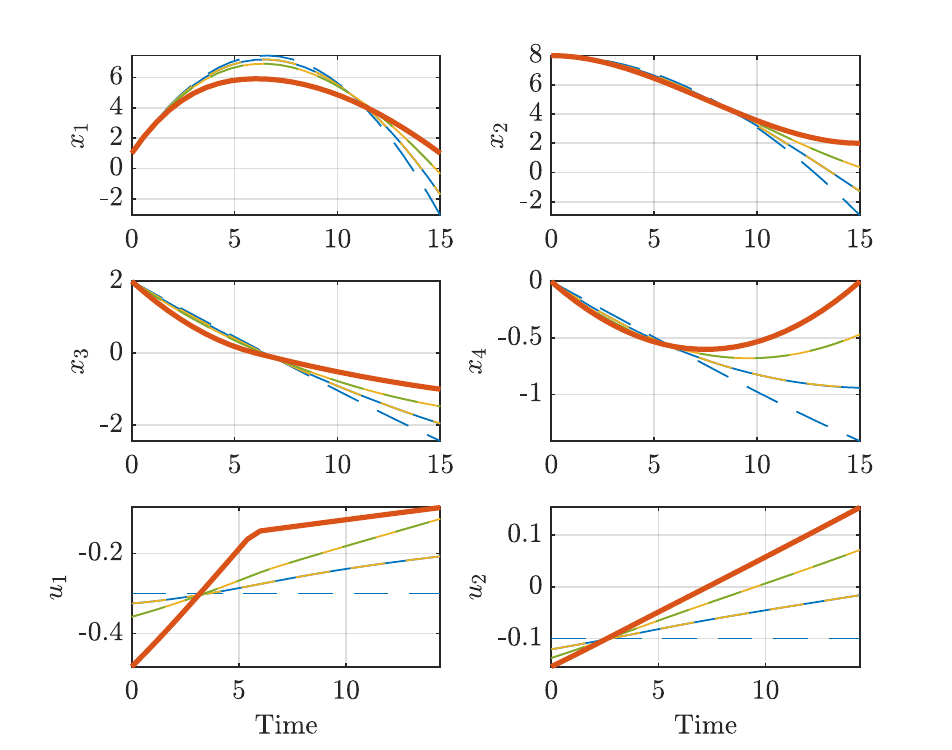}
	\caption{State and control during successive solutions. The dashed lines are the reference $\hat{x}^i$ and $\hat{u}^i$, and the solid lines are the mean state and control after the $i$\textsuperscript{th} step, $\bar{x}^i$ and $\bar{u}^i$. The first iteration is shown in blue and the final iteration is shown in bold.} \label{fig:iterations}
\end{figure}

\begin{figure}
	\centering
	% trim={<left> <lower> <right> <upper>}
	\includegraphics[trim={0.0in 0.1in 0.0in 0.15in}, clip]{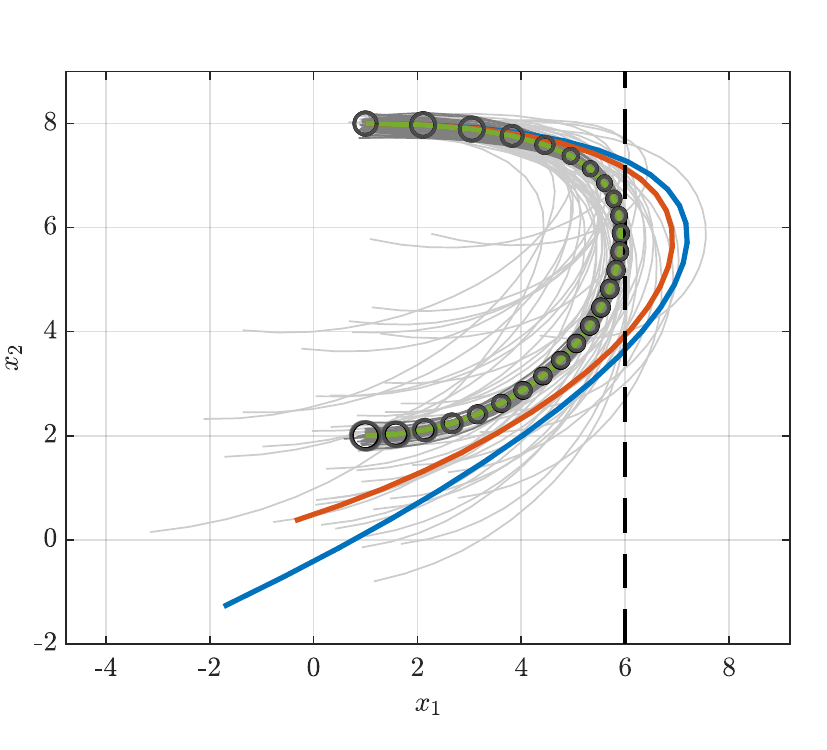}
	\caption{Successive iterations are shown by colored lines, with the initial iteration in blue. The chance constraint is shown by the black dashed line, and 90\% confidence ellipses are shown in black and gray. The black confidence ellipses are computed from the linear analysis and the gray ellipses are computed from Monte Carlo, the dark gray and light gray trajectories are a subset of the Monte Carlo trails for closed and open-loop control, respectively. \label{fig:2d-monte-carlo} }
%	\caption{Successive iterations and Monte Carlo results. Iterations are shown by colored lines, with the initial iteration in blue. The chance constraint is shown by the dashed line, and 90\% confidence ellipses are shown in black and gray. The black confidence ellipses are computed from the linear analysis and the gray ellipses are computed from the Monte Carlo simulations, the dark gray trajectories are a subset of the Monte Carlo trails, and the light gray trajectories are open loop Monte Carlo samples.}
\end{figure}

In this section we apply the iCS algorithm to control a double integrator subject to a quadratic drag force. Let the position $\xi \in \Re ^ 2$ and velocity $v \in \Re ^ 2$ be described by the stochastic system
\begin{align}
	\d \xi_t &= v_t \d t, \\
	\d v_t &= u_t - c_d \norm{v_t} v_t + \gamma \d w_t,
\end{align}
where $c_d > 0$ is the drag coefficient and $\gamma > 0$ is a noise scale parameter. In terms of the state $x = (\xi, v) \in \Re ^ 4$, the dynamics can be written as
\begin{equation}
	\d x_t = [A x_t + B u_t + f_d(x_t)] \d t + G \d w_t,
\end{equation}
where $f_d$ represents the nonlinear drag dynamics and where
\begin{equation}
	A = \begin{bmatrix}
		0_2 & I_2 \\
		0_2 & 0_2
	\end{bmatrix}, \quad B = \begin{bmatrix}
		0_2 \\ I_2
	\end{bmatrix}.
\end{equation}
%The partial derivatives of the drag with respect to the states and controls are
%\begin{equation}
%	\dfrac{\partial f_d}{\partial r} = 0_2 \quad \dfrac{\partial f_d}{\partial v} = - c_d \left(\frac{v v\t}{\norm{v}} + I_2 \norm{v} \right) \quad \dfrac{\partial f_d}{\partial u} = 0_{2}
%\end{equation}
Linearizing about a reference velocity $\hat{v}^i_\tau$, we obtain
%\begin{equation}
%	A^i_\tau = \hat{\sigma}^i \bigg\{ \begin{bmatrix}
%		0_2 & I_2 \\
%		0_2 & 0_2
%	\end{bmatrix} - c_d E_2 \left(\frac{\hat{v}^i_\tau \hat{v}_\tau ^{i \addt} }{\norm{\hat{v}^i_\tau}} + I_2 \norm{\hat{v}^i_\tau} \right) E_2\t \bigg\},
%\end{equation}
\begin{equation}
	A^i_\tau = \sigma \bigg\{ A - c_d E_2 \left(\frac{\hat{v}^i_\tau \hat{v}_\tau ^{i \addt} }{\norm{\hat{v}^i_\tau}} + I_2 \norm{\hat{v}^i_\tau} \right) E_2\t \bigg\},
\end{equation}
where $E_2 \t = \begin{bmatrix}
	0_2, I_2
\end{bmatrix}$, $B^i_\tau = \sigma B$,
%\begin{equation}
%	B^i_\tau = \hat{\sigma}^i \begin{bmatrix}
%		0_2 \\ I_2
%	\end{bmatrix},
%\end{equation}
and $r^i_\tau$ is given as in (\ref{eq:linearized-d-r}). In addition, we enforce the the chance constraint
\begin{equation} \label{eq:example-chance-constraint}
	\Pr (\norm{e_1 \xi_\tau}_1 \leq 6) \geq 1 - 0.1,
\end{equation}
where $e_1 = \begin{bmatrix}
	1, 0
\end{bmatrix}$. The initial state is normally distributed with mean and covariance
\begin{equation}
	\bar{x}_0 = \begin{bmatrix}
		1, 8, 2, 0
	\end{bmatrix} \t, \quad P_{x_0} = 0.01 \times I,
\end{equation}
and the terminal distribution is constrained by the mean and covariance
\begin{equation} \label{eq:example-final-values}
	\bar{x}_f = \begin{bmatrix}
		1, 2, -1, 0
	\end{bmatrix} \t, \quad P_{x_f} = 0.1 \times I.
\end{equation}
We set the drag coefficient $c_d = 0.005$ and the noise scale $\gamma = 0.01$. For the solution, we let the number of discrete steps $N = 25$, terminal mean error weight $w_{x_f} = 1000$, and time scale $\sigma = 15$. The mean cost function was
%$	\ell(x_\tau, u_\tau) = 10 \norm{u_\tau} ^ 2$,
\begin{equation}
	\ell(x_\tau, u_\tau) = 10 \norm{u_\tau} ^ 2,
\end{equation}
the weight matrices were $Q_{x,\tau} \equiv 5 I$ and $Q_{u,\tau} \equiv I$, and the algorithm was seeded with the initial guess
%$ \hat{u}^1_k \equiv \begin{bmatrix}
%        -0.3 & -0.1
%    \end{bmatrix} \t$.
\begin{equation}
    \hat{u}^1_k \equiv \begin{bmatrix}
        -0.3 & -0.1
    \end{bmatrix} \t.
\end{equation}
Since the initial guess violates the chance constraint (\ref{eq:example-chance-constraint}), we relax the chance constraint for the first iteration and tighten it to the final constraint over the first several iterations. The algorithm converged in five iterations, and solutions for each iteration are shown in Figure \ref{fig:iterations}. Samples from a 5,000 trial Monte Carlo simulation are shown in Figure \ref{fig:2d-monte-carlo}. The maximum probability of constraint violation was at step $k = 11$, with 9.14\% of states having $\norm{e_1 \xi_k}_1 \geq 6$, which is below the limit of 10\% set in (\ref{eq:example-chance-constraint}).
%, which can be attributed the effects of the nonlinear term and to the statistical approximation of the Monte Carlo.
Also from the Monte Carlo simulation, the final state mean was
%\begin{equation}
%    \bar{x}_f = \begin{bmatrix}
%        1.000 & 2.004 & -1.000 & 0.001
%    \end{bmatrix},
%\end{equation}
%\begin{equation}
%    P_{x_f} = \begin{bmatrix}
%    0.039  & -0.000  &  0.005 &  -0.000 \\
%    -0.000 &   0.032 &  -0.000  &  0.006 \\
%    0.005  & -0.000  &  0.001 &  -0.000 \\
%    -0.000   & 0.006  & -0.000  &  0.002 \\
%    \end{bmatrix},
%\end{equation}
\begin{equation}
    \bar{x}_f = \begin{bmatrix}
    	1.004 &
    	1.997 &
   		-1.000 &
   		-0.001
    \end{bmatrix} \t ,
\end{equation}
which is very close to the specified value in (\ref{eq:example-final-values}), and the covariance
\begin{equation}
    P_{x_f} = \begin{bmatrix}
    0.018  & -0.001  &  0.004  & 0.000 \\
   -0.001 &  0.016 &  0.000 &   0.004 \\
    0.004  & 0.000  &  0.001  & 0.000 \\
   0.000 &   0.004 &  0.000  &  0.001 \\
    \end{bmatrix}
\end{equation}
is less than the upper bound specified in (\ref{eq:example-final-values}).
%The mean is very close to the prescribed value in (\ref{eq:example-final-values}), and the covariance is less than the maximum value in (\ref{eq:example-final-values}).
%which are very close to the prescribed values in (\ref{eq:example-final-values}).
% \begin{figure}
% 	\centering
% 	\includegraphics[trim={0.1in 0.1in 0.1in 0.1in}, clip, width=3.4in]{mc3x2}
% %		\includegraphics[width=3in]{iterations3x2}
% 	\caption{fo}
% \end{figure}

\section{CONCLUSION} \label{sec:conclusion}

%- overview
%- effect of free final time
%- chance constraints, conservatism
%- initialization
%- future work

In this paper we presented an algorithmic solution to the chance constrained nonlinear CS problem.
We began by approximating the original nonlinear stochastic system by a linear discrete stochastic system, and then we formulated the linear CS problem as a deterministic optimization problem.
Next, in order for the linearized problem formulation to be a reasonable approximation of the original nonlinear problem, we constrained the trajectory at each iteration within a probabilistic trust region about the trajectory from the previous iteration. The size of the trust region depends on the nonlinearity of the dynamics, and therefore convergence properties are problem-specific.

%Since the state covariance and the control covariance are nonconvex in the time dilation coefficient, we approximated the constraints that depend on the covariance by using the value of the time dilation from the previous iteration. As a consequence, these constraints may not be satisfied after some iterations, but a converged solution will satisfy all of the constraints.

Since the proposed iCS algorithm linearizes the dynamics at each iteration, an initial trajectory must be given for the first iteration of the algorithm. At the same time, the difference in the trajectories between iterations is constrained within a trust region, and so a poor initialization may cause the first step to be infeasible. In the numerical example we addressed this problem by relaxing the chance constraints in the first iterations of the algorithm. Since the chance constrained region is assumed to be a convex polytope, the region can be easily expanded by scaling the inequality constraints. Another solution would be to first solve a deterministic optimization problem with tightened inequality constraints representing a worst-case chance constrained region. In this case, the iCS algorithm would be used to improve the solution from the deterministic problem by adding the closed-loop system statistics to the optimization. The latter approach could be applied to problems such as planetary entry and powered descent by iterating on a given reference trajectory that is to be tracked in the presence of uncertainty.
% Note: cut out reference Cianciolo2017

In future work we plan to apply iCS to problems in entry, descent, and landing (EDL) with nonlinear dynamics, such as entry and powered descent \cite{Ridderhof2019}. 
Another extension to this work would be the addition of time-varying chance constraints that are satisfied for all time, rather than for each time, while not being overly conservative.

%%%%%%%%%%%%%%%%%%%%%%%%%%%%%%%%%%%%%%%%%%%%%%%%%%%%%%%%%%%%%%%%%%%%%%%%%%%%%%%%
%\section{ACKNOWLEDGMENTS}

%This work was supported by a NASA Space Technology Research Fellowship.

%%%%%%%%%%%%%%%%%%%%%%%%%%%%%%%%%%%%%%%%%%%%%%%%%%%%%%%%%%%%%%%%%%%%%%%%%%%%%%%%

\nocite{*}
\bibliography{cdc19}
\bibliographystyle{ieeetr}

%\begin{thebibliography}{99}

%\end{thebibliography}

\end{document}